\documentclass[11pt]{article}
\usepackage{amssymb,amsmath,latexsym}

\usepackage{color}

\usepackage{epsfig}

\newtheorem{thm}{Theorem}[section]
\newtheorem{defn}[thm]{Definition}

\newtheorem{lemma}[thm]{Lemma}

\newtheorem{remark}[thm]{Remark}

\newtheorem{assumption}[thm]{Assumption}

\newcommand{\pf}{\noindent{\bf Proof.} }
\def\qed{{\hfill $\Box$ \bigskip}}

\newcommand\cbrk{\text{$]$\kern-.15em$]$}}
\newcommand\opar{\text{\,\raise.2ex\hbox{${\scriptstyle
|}$}\kern-.34em$($}}
\newcommand\cpar{\text{$)$\kern-.34em\raise.2ex\hbox{${\scriptstyle |}$}}\,}

\def\wh{\widehat}
\def\wt{\widetilde}
\def\<{\langle}
\def\>{\rangle}

\def\E{{\mathbb E}}

\newcommand\bL{\mathbb{L}}
\newcommand\bR{\mathbb{R}}
\newcommand\bH{\mathbb{H}}

\newcommand\bD{\mathbb{D}}
\newcommand\bS{\mathbb{S}}

\newcommand\cA{\mathcal{A}}
\newcommand\cB{\mathcal{B}}

\newcommand\cF{\mathcal{F}}
\newcommand\cH{\mathcal{H}}

\newcommand\cP{\mathcal{P}}

\newcommand\cM{\mathcal{M}}

\def\wh{\widehat}
\def\wt{\widetilde}

\newcommand{\mysection}[1]{\section{#1}
\setcounter{equation}{0}}

\begin{document}

\title{\bf  An $L^p$-theory of
non-divergence form SPDEs
driven by L\'e{}vy processes}

\author{Zhen-Qing Chen \qquad \hbox{\rm and} \qquad Kyeong-Hun Kim}

\date{}

\author{Zhen-Qing Chen\footnote{Department of Mathematics, University of Washington,
Seattle, WA 98195, USA, \,\, zchen@math.washington.edu. The research
of this author is supported in part by NSF Grant DMS-0906743.}
\qquad \hbox{\rm and} \qquad Kyeong-Hun Kim\footnote{Department of
Mathematics, Korea University, 1 Anam-dong, Sungbuk-gu, Seoul, South
Korea 136-701, \,\, kyeonghun@korea.ac.kr. The research of this
author was supported by Basic Science Research Program through the
National Research Foundation of Korea(NRF) funded by the Ministry of
Education, Science and Technology (20090087117)}}

\maketitle

\begin{abstract}

In this paper we  present  an $L^p$-theory for
 the stochastic partial differential equations (SPDEs in abbreciation) driven by L\'e{}vy processes.  Existence and
uniqueness of solutions in  Sobolev spaces are obtained.  The
coefficients of SPDEs under consideration
 are  random functions depending on
time and space variables.

\vspace*{.125in}


\vspace*{.125in}


\end{abstract}

\mysection{Introduction}

Let $(\Omega,\cF,P)$ be a complete probability space,
$\{\cF_{t},t\geq0\}$ be an increasing filtration of $\sigma$-fields
$\cF_{t}\subset\cF$, each of which contains all $(\cF,P)$-null sets.
We assume that on $\Omega$ we are given independent  one-dimensional L\'evy processes
 $Z^{1}_{t},Z^{2}_{t},...$
relative to $\{\cF_{t},t\geq0\}$. Let $\cP$ be the predictable
$\sigma$-field generated by $\{\cF_{t},t\geq0\}$.

 In this article we are dealing with
 $W^{n,p}$-theory
 of the stochastic   partial differential equation
\begin{equation}   \label{e:1.1}
du=(a^{ij}u_{x^ix^j}+b^iu_{x^i}+cu+f)dt
+(\sigma^{ik}u_{x^i}+\mu^{k}u+g^k)dZ^k_t
\end{equation}
given for $t\geq0$  and $x\in \bR^d$.  Here $p\in [2,\infty)$ and $n\in \bR$.  Indices $i$ and $j$
go from $1$ to $d$, and $k$ runs through $\{1,2,...\}$ with the
summation convention on $i,j,k$ being enforced. The coefficients
$a^{ij}$, $b^i$, $c, \sigma^{ik}, \mu^{k}$ and the free terms $f, g^k$ are
random functions depending on $(t,x)$.

Demand for a general theory of stochastic partial differential equations (SPDEs) driven by
L\'evy processes is obvious
 when
we model  a natural
phenomenon with randomness and jumps.
The main objective of this paper is to establish
unique solvability in  Sobolev spaces for SPDEs \eqref{e:1.1}.

If $\{Z^k, k\geq 1\}$
are independent one-dimensional Wiener processes,
$L^p$-theory  for  SPDE \eqref{e:1.1}
 has been well studied. An $L^p$-theory of SPDEs with Wiener processes defined on $\bR^n$ was first introduced by Krylov in \cite{Kr99}.
  Subsequently,  Krylov and Lototsky \cite{KL1, KL2}
  developed an $L^p$-theory of such equations
in half space $\bR^n_+$ with constant coefficients.
These results were later extended to SPDEs with variable coefficients defined in bounded domains of $\bR^n$ by several
authors, see, for instance, \cite{KK, Kim03, Lo}.

However very little is known when $\{Z^k, k\geq 1\}$
are general discontinuous L\'evy processes. As
  far as we know,  most previous
   work on SPDEs driven by L\'evy processes
   deal with equations with non-random coefficients independent of $t$ and, moreover, $\sigma^{ik}$ have always been assumed to be zero, consequently first derivatives of solutions were not allowed to appear in the stochastic part. More precisely, the typical type of equations appearing in the previous works (see
   \cite{AWZ, F, M, RZ} and references therein) is of the following type:
\begin{equation}  \label{e:1.2}
du=(Au+f)dt +\sum_{k=1}^n g^k(u)dZ^k_t,
\end{equation}
where $A(t)$ is the generator of certain
 semigroup.
and the function $g^k$ satisfies certain continuity conditions.

Our approaches are different from those in \cite{AWZ, F, M, RZ}. We
 adopt  analytic approaches introduced  by
Krylov in \cite{Kr99}. Our results are new even for the stochastic
heat equation
$$
du=\Delta u \, dt +g \, dZ_t,
$$
since we
 establish
  unique solvability result in $L^p(\Omega\times [0,T],H^{n}_p)$
 for every $p\in [2,\infty)$ and $n\in \bR$, not just
 in $L^2(\Omega\times [0,T],H^1_2)$ space.
(The definition of Sobolev space $H^n_p=W^{n,p}$
 will be given in next section.)
  This allows us to
obtain  various regularity results of solutions. See Remark \ref{main remark}(ii).

Our main results are stated in section 2 and consist of
Theorem \ref{main thm 1} ($L^2$-theory) and Theorem \ref{main thm 2} ($L^p$-theory, $p>2$). In section 3 we deal with equations with constant coefficients,
and in section 4 we prove Theorem \ref{main thm 1} and Theorem \ref{main thm 2}.

We end the introduction with some notation.
As usual $\bR^{d}$
stands for the Euclidean space of points $x=(x^{1},...,x^{d})$,
 and
$B_{r}(x)=\{y\in\bR^{d}:|x-y|<r\}$.
 For $i=1,...,d$, multi-indices $\alpha=(\alpha_{1},...,\alpha_{d})$
 with  $\alpha_{i}\in\{0,1,2,...\}$
 and functions $u(x)$ on $\bR^d$, we set
$$
 D_{i}u:=u_{x^{i}}:=\partial u/\partial x^{i}, \qquad
D^{\alpha}u:=D_{1}^{\alpha_{1}}\cdot...\cdot D^{\alpha_{d}}_{d}u,
\qquad |\alpha|:=\alpha_{1}+...+\alpha_{d}.
$$
For $a, b\in \bR^d$, we define $a\wedge b:=\min \{a, \, b\}$
and $a\vee b:=\max \{a, \, b\}$.
For $p\geq 1$, we will use $\| u\|_p$ to denote the $L^p$-norm
of $u$ in $ L^p(\bR^d; dx)$. For scalar functions $f, g$ on $\bR^d$,
$(f, g):= \int_{\bR^d} f(x) g(x) dx$.

\mysection{Main results}

For $t\geq 0$ and $A\in \cB(\bR\setminus \{0\})$,  define
$$
N^k(t,A)
:=\# \left\{0\leq s\leq t; \, Z^k_s-Z^k_{s-} \in A \right\}, \quad
\wt {N}^k(t,A):=N^k(t,A)-t\nu_k(A)
$$
where $\nu_k(A):=\E [N_k(1,A)]$ is the L\'e{}vy measure of $Z^k$. By
L\'e{}vy-It\^o decomposition,  there exist constants $\alpha^k,
\beta^k$ and Brownian motion $B^k$ so that
\begin{equation}
                                           \label{eqn 2.27}
Z^k(t)=\alpha^kt +\beta^k B^k_t+\int_{|z|<1}z \wt
{N}^k(t,dz)+\int_{|z|\geq 1} z N^k(t, dz).
\end{equation}
For simplicity, throughout this paper,
we assume $\alpha^k=\beta^k=0$. Consider the
following equation for random function $u(t, x)$
on $\Omega \times [0, \infty)\times \bR^d$:
\begin{equation}      \label{e:2.1}
du= \left( a^{ij}u_{x^ix^j}+b^iu_{x^i}+ cu+ f \right) dt + \left(
\sigma^{ik} u_{x^i}+\mu^{k}u+g^{k} \right) dZ^{k}_{t}
\end{equation}
in the weak sense.
Precise definition of weak solution to \eqref{e:2.1} will be
given in Definition \ref{D:2.5} below.
 Here    $i$ and $j$ go from $1$ to $d$, and $k$
runs through $\{1,2,\cdots\}$. The coefficients $a^{ij},
b^i,c,\sigma^{ik},\mu^k$ and the free terms $f,g^k$ are random
functions depending on $t>0$ and $x\in \bR^d$.

\begin{assumption}        \label{assmp 3}
 $p\in [2,\infty)$ and for each $k$,
$$
\wh {c}_{k,p}:=\left(\int_{\bR} |z|^p \nu_k(dz)\right)^{1/p}<\infty.
$$
\end{assumption}
For $n=0,1,2,...$, define Sobolev space
$$
H^n_p:=H^n_p(\bR^d) =\left\{u: u, Du,...,D^n u\in L^p (\bR^d) \right\},
$$
where $D^ku$ are derivatives in the distributional sense.
 In literature, $H^n_p$ is also denoted as
$W^{n,p}(\bR^d)$. In general, for $\gamma \in \bR$ define the space
$H^{\gamma}_p=H^{\gamma}_p(\bR^d)=(1-\Delta)^{-\gamma/2}L^p$ (called
the space of Bessel potentials or the Sobolev space with fractional
derivatives) as the set of all distributions $u$ such that
$(1-\Delta)^{\gamma/2}u\in L^p$. For $u\in H^\gamma_p$, we define
\begin{equation}   \label{e:2.2}
\|u\|_{H^{\gamma}_p}:=\|(1-\Delta)^{\gamma/2}u\|_p
:=\|\cF^{-1}[(1+|\xi|^2)^{\gamma/2}\cF(u)(\xi)]\|_p,
\end{equation}
where $\cF$ is the Fourier transform.
 Let $\cP^{dP\times dt}$ be the completion of $\cP$ with respect to  $dP \times dt$.
 Denote by $\bH^{\gamma}_p(T)$ the space
  of all $\cP^{dP\times dt}$-measurable processes
$u:[0,T]\times \Omega\to H^{\gamma}_p$ so that
$$
\|u\|_{\bH^{\gamma}_p(T)}:= \left( \E \left[ \int^T_0
\,\|u\|^p_{H^{\gamma}_p}\,dt \right] \right)^{1/p}<\infty.
$$
 For fixed $p\geq 2$, define
$$
\wh{c}_{k}:=\wh{c}_{k,2} \vee \wh{c}_{k,p}.
$$
Note that $\wh c_k=\wt c_{k, 2}$ when $p=2$, and
for   $2< q<p$, by H\"older's inequality,
$$ \wh c_{k, q} \leq \left(\int_{\bR} |z|^2 \nu_k(dz)\right)^{(p-q)/(q(p-2))}
\left(\int_{\bR} |z|^p \nu_k(dz)\right)^{(q-2)/(q(p-2))}
\leq \wh c_k.
$$
For $\ell^2$-valued processes $g=(g^1,g^2,...)$,
 we say $g\in
\bH^{\gamma}_p(T,\ell^2)$ if $g^k\in \bH^{\gamma}_p(T)$ for every
$k\geq 1$ and
\begin{equation}   \label{e:2.3}
\|g\|_{\bH^{\gamma}_p(T,\ell^2)}:= \left( \E \left[\int^T_0\|\,
|(1-\Delta)^{\gamma/2}\wh{g}|_{\ell^2} \,
 \|^p_p\,dt \right] \right)^{1/p}<\infty,
\end{equation}
where
$\wh g= (\wh g_1, \wh g_2, \wh g_3, \cdots ):=
(\wh{c}_{1}g^1,\wh{c}_{2}g^2,\wh{c}_{3}g^3,\cdots)$.
Finally, we say $u_0\in U^{\gamma}_p$ if $u_0$ is
$\cF_0$-measurable and
$$
 \|u_0\|_{U^{\gamma}_p}:=\left( \E \left[
\|u_0\|^p_{H^{\gamma-(2/p)}_p} \right]\right)^{1/p}<\infty.
$$

\begin{remark}   \label{remark 11.22} \rm
\rm It follows from \eqref{e:2.2} that for any $\mu, \gamma\in \bR$,
the operator $(1-\Delta)^{\mu/2}:H^{\gamma}_p\to H^{\gamma-\mu}_p$
is an isometry. Indeed,
$$
\|(1-\Delta)^{\mu/2}u\|_{H^{\gamma-\mu}_p}
=\|(1-\Delta)^{(\gamma-\mu)/2}(1-\Delta)^{\mu/2}u\|_p
=\|(1-\Delta)^{\gamma/2}u\|_p=\|u\|_{H^{\gamma}_p}.
$$
\end{remark}

\begin{remark}\label{R:2.4}\rm \begin{description}
\item{(i)} Let $\cM^{\gamma}_p(T)$ denote the set of all $H^{\gamma}_p$-valued $\{\cF_t\}$-adapted
processes $u(t)$ that are $\cF \otimes \cB(0,T)$-measurable and
satisfy
$$
\E \left[ \int^T_0\|u\|^p_{H^{\gamma}_p}\, dt \right] <\infty.
$$
Then by Theorem 2.8.2 in \cite{kr}, $ \cM^{\gamma}_p(T) \subset
\bH^{\gamma}_p(T)$.

\item{(ii)} Note that
under Assumption \ref{assmp 3} for $p=2$,
  $Z^k  $ is a square integrable martingale for each $k\geq 1$. For every
  $\cP^{dP\times dt}$-measurable
 process $H\in L^2(\Omega\times [0,T])$, $M_t:=\int_0^t H_s dZ^k_s$ is a square integrable
martingale with
\begin{equation}\label{e:2.4}
\E [M_t^2] = \E \left[ \int_0^t H_s^2 \, d[Z^k]_s\right]
 =  \wh{c}_{k,2}^2  \E \left[ \int_0^t H_s^2 ds\right].
 \end{equation}
So for any  $g\in \bH^{\gamma}_p(T,\ell^2)$ and $\phi\in C^\infty_c
(\bR^d)$,
 \begin{eqnarray*}
&& \sum_{k=1}^{\infty}\wh{c}_k^2 \, \E \left[ \int^T_0(g^k,\phi)^2ds \right]\\
&=&\sum_{k=1}^{\infty} \E \left[ \int_0^T
((1-\Delta)^{\gamma/2}\wh{c}_kg^k,
(1-\Delta)^{-\gamma/2}\phi)^2 \, ds \right] \\
  &\leq&
\|(1-\Delta)^{-\gamma/2}\phi\|_1 \, \E \left[ \int^T_0 \Big (
\sum_k |(1-\Delta)^{\gamma/2}\wh{c}_kg^k|^2, \,
|(1-\Delta)^{-\gamma/2}\phi| \Big) \, ds \right] \\
&\leq& \|(1-\Delta)^{-\gamma/2}\phi\|_1 \,
\|(1-\Delta)^{-\gamma/2}\phi\|_q \, \E \left[ \int^T_0 \Big\|\,
 \sum_k |(1-\Delta)^{\gamma/2}\wh{c}_kg^k|^2    \Big\|_{p/2}\,
 ds \right]\\
&\leq&  \|(1-\Delta)^{-\gamma/2}\phi\|_1 \,
\|(1-\Delta)^{-\gamma/2}\phi\|_q \, T^{1-\frac2{p}} \, \|
g\|_{\bH^\gamma_p(T, l^2)}^2 <\infty ,
\end{eqnarray*}
where $q=p/(p-2)$.
 Thus in view of \eqref{e:2.4},
 the series of stochastic integral
$\sum_{k=1}^{\infty}\int^t_0(g^k,\phi)dZ^k_s$ defines a square
integrable martingale on $[0, T]$, which is right continuous with
left limits.
\end{description}
\end{remark}

\begin{defn}\label{D:2.5}
Write $u \in \cH^{\gamma+2}_p(T)$ if $u\in \bH^{\gamma+2}_p(T)$ with
$u(0)\in U^{\gamma}_p$,  and for some $f\in \bH^{\gamma}_p(T)$ and
 $ g\in \bH^{\gamma+1}_p(T,\ell^2)$
$$
du=f dt +g^k dZ^{k}_t, \quad \hbox{for }  t\in [0,  T]
$$
 in the distributional sense, that is, for any $\phi\in C^\infty_c (\bR^d)$,
\begin{equation}
                \label{eqn 11.16}
(u(t),\phi)=(u(0),\phi)+\int^t_0(f,\phi)dt +
\sum_k\int^t_0(g^k,\phi)dZ^{k}_t
\end{equation}
holds for all $t\leq T$ $a.s.$.
Define
$$
\bD u :=f, \quad \bS u :=g,
$$
and define
$$
\|u\|_{\cH^{\gamma+2}_p(T)}:=\|u\|_{\bH^{\gamma+2}_p(T)}
 +\|\bD u\|_{\bH^{\gamma}_p(T)} +
 \|\bS u\|_{\bH^{\gamma+1}_p(T,\ell^2)}
+\|u(0)\|_{U^{\gamma+2}_p}.
$$
\end{defn}

\begin{thm}
                        \label{theorem banach}
For any $p\in [2, \infty)$,$\gamma\in \bR$ and $T>0$,
$\cH^{\gamma+2}_p(T)$ is a Banach space with norm $\|
\cdot\|_{\cH^{\gamma+2}_p(T)}$. Moreover, there is a constant
$c=c(d,p)>0$, independent of $T$,
 such that for every $u\in
\cH^{\gamma+2}_p(T)$,
\begin{equation}
                                         \label{e:2.6}
\E \left[ \sup_{t\leq T}\|u(t)\|^p_{H^{\gamma}_p} \right] \leq c
\left( \|\bD u\|^p_{\bH^{\gamma}_p(T)}+
 \|\bS u\|^p_{\bH^{\gamma }_p(T,\ell^2)}
+\E \left[ \|u_0\|^p_{H^{\gamma}_p} \right]\right) .
\end{equation}
Consequently, for each  $t>0$,
\begin{equation}
                    \label{new 2}
\|u\|^p_{\bH^{\gamma}_p(t)}\leq
c\int^t_0\|u\|^p_{\cH^{\gamma+2}_p(s)} ds.
\end{equation}
\end{thm}

\pf In view of  Remark \ref{remark 11.22} it suffices to prove the
theorem  for $\gamma=0$.  First we prove (\ref{e:2.6}). Let
$du=fdt +g^kdZ^k_t$ with $u(0)=u_0$.  Assume that $g^k=0$ for all
$k\geq N_0$ and  $g^k$ is of the type
\begin{equation}
                               \label{e:12.5}
g^k(t,x)=\sum_{i=0}^m I_{(\tau^{k}_{i},\tau^{k}_{i+1}]}(t)g^{ki}(x),
\end{equation}
where $\tau^k_i$ are bounded stopping times and $g^{ki}\in
C^{\infty}_c(\bR^d)$. Define
$$
v(t,x)=\sum_{k=1}^{\infty} \int^t_0 g^kdZ^k_s.
$$
Then by Burkholder-Davis-Gundy inequality (used twice)
and monotone convergence
theorem,
\begin{eqnarray*}
&&\E \left[ \sup_{s\leq t}|v(s,x)|^p \right] \\
&\leq&  c \E
\left[\left(\sum_{k=1}^{\infty}\int^t_0\int|g^k(s,x)|^2|z|^2
N^k(ds,dz)\right)^{p/2}\right]\\
&=&c\lim_{N\to \infty}\E
\left[ \left(\sum_{k=1}^{\infty}\int^t_0\int_{|z|\leq N}|g^k(s,x)|^2|z|^2N^k(ds,dz)\right)^{p/2}\right]\\
&\leq& c \lim_{N\to \infty}\E \left[
\left(\sum_{k=1}^{\infty}\int^t_0\int_{|z|\leq N}|g^k(s,x)|^2|z|^2
\wt{N}^k(ds,dz)\right)^{p/2}\right] \\
&& +\, c\E\left[ \left( \sum_{k=1}^{\infty}\int^t_0\int_{\bR}
\sum_{k=1}^\infty |g^k(s, x)|^2\,|z|^2\nu_k(dz)ds\right)^{p/2}
\right] \\
&\leq& c\lim_{N\to \infty}\E \left[\left(
\sum_{k=1}^{\infty}\int^t_0\int_{|z|\leq N}
 |g^k(s, x)|^4|z|^4 N^k(ds,dz)\right)^{p/4}\right]
+c\E\left[\left(\int^t_0\sum_{k=1}^{\infty}|\wh{g}^k (s, x)|^2 \, ds
\right)^{p/2}\right].
\end{eqnarray*}
Recall that for any $q>1$, $(\sum |a_n|^q)^{1/q}\leq \sum |a_n|$.
Thus if $2< p \leq 4$, then
\begin{eqnarray*}
&&  \E \left[\left(\sum_{k=1}^{\infty}\int^t_0\int_{|z|\leq N}|g^k(s, x) |^4 \, |z|^4 N^k(ds,dz)\right)^{p/4}\right] \\
&\leq &\E\left[\left( \sum_k \sum_{0\leq s\leq
t}|g^k (s,x)|^4 \, |\Delta Z^k_s|^4\right)^{p/4}\right]
  \leq  \E \left[ \sum_k\sum_{0\leq
s\leq t}|g^k (s,x)|^p \, |\Delta Z^k_s|^p \right] \\
&=&\E \left[ \int^t_0\sum_{k=1}^{\infty}|\wh{c}_{k,p}g^k (s, x)|^p
 \, ds \right] \leq
\E \left[ \int^t_0\left(\sum_{k=1}^{\infty}|\wh{c}_{k,p}g^k (s, x)|^2\right)^{p/2}\, ds \right]
\end{eqnarray*}
If $4<p\leq 8$ then
\begin{eqnarray*}
&& \E \left[ \left(\sum_{k=1}^{\infty}\int^t_0\int_{|z|\leq N}|g^k(s, x)|^4 \, |z|^4 N^k(ds,dz)\right)^{p/4}\right]\\
&\leq&\E\left[ \left(\sum_{k=1}^{\infty}\int^t_0\int_{|z|\leq
N}|g^k(s,x)|^4 \, |z|^4
\wt{N}^k(ds,dz)+\sum_{k=1}^{\infty}\int^t_0\int_{|z|\leq N}|g(s,x)|^4 \, |z|^4\nu_k(dz)ds\right)^{p/4}\right] \\
&\leq& c\E \left[\left(\sum_{k=1}^{\infty}\int^t_0\int_{|z|\leq
N}|g^k (s, x)|^8|z|^8 N^k(ds,dz)\right)^{p/8} +\left(
\int^t_0\sum_{k=1}^{\infty}|\wh{g}^k (s, x)|^4
ds \right)^{p/4}\right]\\
&\leq& c\E\left[ \int^t_0\sum_{k=1}^{\infty}|\wh{c}_{k,p}g^k (s,
x)|^p \, ds+ \left(\int^t_0\sum_{k=1}^{\infty}|\wh{g}^k (s, x)|^4 ds
\right)^{p/4}\right].
\end{eqnarray*}
 Similarly, in general, for $p\in (2^{n-1}, 2^n]$,
 \begin{eqnarray*}
  &&  \E\left[\left( \sum_{k=1}^{\infty}\int^t_0\int|g^k (s, x)|^2\, |z|^2 N^k(dz,ds)\right)^{p/2}\right] \\
 & \leq & c  \, \E \left[ \sum_{j=1}^n \left(\int^t_0\sum_k|\wh{g}^k (s, x)|^{2j} ds\right)^{p2^{-j}} \right]  +c\, \E \left[\int^t_0\sum_k|\wh{g}^k(s, x)|^pds\right].
 \end{eqnarray*}
 Also since for each $2\leq q\leq
 p$,
$$
\left(\int^t_0\sum_{k= 1}^\infty |\wh{g}^k (s, x)|^q ds\right)^{1/q}\leq c(q)\left(\left(\int^t_0\sum_{k=1}^\infty
 |\wh{g}^k(s, x)|^2 ds\right)^{1/2}+\left(\int^t_0
 \sum_{k=1}^\infty |\wh{g}^k(s,x)|^p ds\right)^{1/p}\right),
$$
we get
\begin{eqnarray}
                             \label{e: 12.6.4}
&&  \E\left[\left(\sum_{k=1}^{\infty}\int^t_0\int|g^k(s,x)|^2\,|z|^2
N^k(dz,ds)\right)^{p/2}\right] \nonumber \\
 &\leq&
c(p) \, \E \left[\left(\int^t_0\sum_k|\wh{g}^k(s, x)|^2ds\right)^{p/2}+ \int^t_0\sum_k|\wh{g}^k(s, x)|^p\, ds \right]
\end{eqnarray}
and
\begin{equation}
                           \label{e:12.6.2}
 \E \left[ \sup_{s\leq t}|v(s,x)|^p \right]\leq
 c(p)\, \E\left[ \left(\int^t_0\sum_{k=1}^\infty |\wh{g}^k (s, x)|^2ds \right)^{p/2}+ \int^t_0\sum_k|\wh{g}^k(s, x)|^pds\right] .
 \end{equation}
By integrating over $\bR^d$, we get
\begin{equation}
                     \label{eqn 11.22.1}
\E \left[\sup_{s\leq t}\|v\|^p_p \right]\leq c(p)
\|g\|^p_{\bH^0_p(t,\ell^2)}.
\end{equation}
Next we show that \eqref{eqn 11.22.1} holds for general
$g\in \bH^0_p(T,\ell^2)$.
 Take a
sequence $g_n\in \bH^0_p(T,\ell^2)$ so that for each fixed $n$,
$g^k_n=0$ for all large $k$ and each $g^k_n$ is of of the type
(\ref{e:12.5}), and $g_n \to g$ in $\bH^0_p(T, \ell^2)$ as $n\to
\infty$. Define $v_n(t,x)=\sum_{k}\int^t_0 g^k_n dZ^k_t$, then
$$
\E \left[ \sup_{s\leq t}\|v_n\|^p_p\right] \leq c(p)
\|g_n\|^p_{\bH^0_p(t,\ell^2)}, \quad \E \left[\sup_{s\leq
t}\|v_m-v_n\|^p_p\right] \leq c(p)
\|g_m-g_n\|^p_{\bH^0_p(t,\ell^2)}.
$$
Thus
 \eqref{eqn 11.22.1} follows by taking $n\to \infty$.
Now note that
$$
d(u-v)=f dt  \quad \hbox{ with} \quad (u-v)(0)=u_0.
$$
Thus it is easy to check that
$$
\E \left[ \sup_{s\leq t}\|u-v\|^p_p \right] \leq c\E \left[
\|u_0\|^p_p\right] +c\E \left[ \int^t_0\|f(s, \cdot) \|^p_p \,
ds\right].
$$
Consequently,
$$
\E \left[ \sup_{s\leq t}\|u\|^p_p \right] \leq
 c
\|f\|^p_{\bH^0_p(t)}+c\|g\|^p_{\bH^0_p(t,\ell^2)}+c\E\|u_0\|^p_{L_p}.
$$
The completeness of the space $\cH^2_p(T)$ easily follows from
\eqref{e:2.6}. The theorem is proved. \qed

Now we introduce the space of point-wise multipliers in
$H^{\gamma}_p$. Fix $\kappa_0>0$. For $r\geq 0$,  define
$r_+=r$ if $r=0,1,2,\cdots$, and $r_+=r+\kappa_0$ otherwise. Also denote $r^+=r+\kappa_0$.  Define
 \begin{equation}\label{e:2.11}
 B^{r}=\begin{cases} B(\bR^d)  \qquad & \hbox{if }  r=0,\\
 C^{r-1,1}(\bR^d) &\hbox{if }  r=1,2, \cdots ,\\
C^{r}(\bR^d) &  \hbox{otherwise},
 \end{cases}
\end{equation}
where $B(\bR^d)$ is the space of
bounded Borel measurable functions
on $\bR^d$, $C^{r-1,1}(\bR^d)$ is the space of $r-1$ times
continuously differentiable functions whose $(r-1)$st order derivatives  are Lipschitz continuous, and
$C^{r}(\bR^d)$ is the usual H\"older space. We also
use the Banach space $B^{r}$ for $\ell^2$-valued functions.
For instance,
 if $g=(g^1,g^2,...)$,  then $|g|_{B^0}=\sup_x |g(x)|_{\ell^2}$ and
$$
|g|_{C^{n-1,1}}=\sum_{|\alpha|\leq n-1} |D^{\alpha}g
|_{B^0}+\sum_{|\alpha|=n-1}\sup_{x\neq y}\frac{
|D^{\alpha}g(x)-D^{\alpha}g(y)|_{\ell^2}}{|x-y|}.
$$

\begin{assumption}
                         \label{assumption 1}
(i) The coefficients $a^{ij},  b^i, c, \sigma^{ik},\mu^k$ are
$\cP\otimes \cB(\bR^d)$-measurable functions.

(ii) $a^{ij}=a^{ji}$, and  the functions $a^{ij}$ and $\sigma^i$ are
uniformly
 continuous in $x$. In other words, for any
 $\varepsilon >0$, there exists $\delta>0$ such that
 whenever  $|x-y|<\delta$,
 $$
 |a^{ij}(t,x)-a^{ij}(t,y)|+|\sigma^i(t,x)-\sigma^i(t,y)
 |_{\ell^2}< \varepsilon .
 $$

(iii) There exist constants $\delta, K>0$ so that
$$
|a^{ij}|+|b^i|+|c|+ |\sigma^{i}|_{\ell^2}+|\mu |_{\ell^2}\leq K,
$$
\begin{equation}
                                      \label{eqn 8.5.1}
 \delta I_{d\times d}\leq (a^{ij}-\alpha^{ij} )\leq (a^{ij})\leq K I_{d\times d},
\end{equation}
where $\alpha^{ij}:= \frac{1}{2} \sum_{k=1}^\infty \wh {c}_{k,2}^2
\sigma^{ik}\sigma^{jk}$ and $I_{d\times d}$ denotes the $(d\times
d)$-identity matrix.
\end{assumption}

Here are  main results of this article. We formulate them into two theorems since our assumptions are
 stronger when $p\neq 2$.

\begin{thm}
                      \label{main thm 1}
Let $ \gamma\in \bR, T>0$ and Assumption \ref{assumption 1} hold.
Also assume  there is a constant $L>0$ so that for each $\omega,t,$
$$
 |a^{ij}(t,\cdot) |_{B^{|\gamma|_+}} +
 |b^i(t,\cdot) |_{B^{|\gamma|_+}}+
 |c(t,\cdot) |_{B^{|\gamma|_+ }}
+ |\sigma^i(t,\cdot)|_{B^{|\gamma+1|_+}}
+|\mu(t,\cdot)|_{B^{|\gamma+1|_+}} \leq L.
$$
Then  for any $f\in \bH^{\gamma}_2(T)$, $g\in \bH^{\gamma+1}_2(T,
\ell^2)$ and $u_0\in U^{\gamma+2}_2$  equation (\ref{e:2.1}) has a
unique solution $u\in \cH^{\gamma+2}_2(T)$, and
\begin{equation}
                            \label{eqn main L_2}
\|u\|_{\cH^{\gamma+2}_2(T)}\leq
c\left(\|f\|_{\bH^{\gamma}_2(T)}+\|g\|_{\bH^{\gamma+1}_2(T,
\ell^2)}+\|u_0\|_{U^{\gamma+2}_2}\right),
\end{equation}
where $c=c(\delta,K,L,\gamma,T)$.
\end{thm}

\begin{remark}
Condition (\ref{eqn 8.5.1}) naturally appears when one writes
It\^o's formula for $|u|^2$, where $u$ is a solution of
(\ref{e:1.1}) (see Lemma 2.8 in \cite{C-K}). Remember we assumed
$\beta^k=0$ in (\ref{eqn 2.27}). If $\beta^k\neq 0$ we need to
replace (\ref{eqn 8.5.1}) by
$$
\delta I_{d\times d}\leq (a^{ij}-\bar{\alpha}^{ij} )\leq
(a^{ij})\leq K I_{d\times d},
$$
where $\bar{\alpha}^{ij}:= \frac{1}{2} \sum_{k=1}^\infty
((\beta^k)^2+\wh {c}_{k,2}^2) \sigma^{ik}\sigma^{jk}$.
\end{remark}

\begin{thm}
                      \label{main thm 2}
Let $p\in (2,\infty), \gamma\in \bR$  and $\varepsilon>0$ be fixed.
Assume Assumption \ref{assumption 1} holds, $\sigma^{ik}=0$ and
 there is a constant $L>0$ so that for each $\omega,t$,
$$
|a^{ij}(t,\cdot) |_{B^{|\gamma|_+}} +
 |b^i(t,\cdot) |_{B^{|\gamma|_+}}+
 |c(t,\cdot) |_{B^{|\gamma|_+ }}
+|\mu(t,\cdot) |_{B^{|\gamma+1|^+}} \leq L.
$$
Then  for any $f\in \bH^{\gamma}_p(T)$, $g\in
\bH^{\gamma+1+\varepsilon}_p(T, \ell^2)$ and $u_0\in
U^{\gamma+2}_p$, equation (\ref{e:2.1}) has a unique solution $u\in
\cH^{\gamma+2}_p(T)$, and
\begin{equation}
                              \label{eqn main L_p}
\|u\|_{\cH^{\gamma+2}_p(T)}\leq
c\left(\|f\|_{\bH^{\gamma}_p(T)}+\|g\|_{\bH^{\gamma+1+\varepsilon}_p(T,
\ell^2)}+\|u_0\|_{U^{\gamma+2}_p}\right),
\end{equation}
where $c=c(\delta,K,L,p,\gamma,T)$.
\end{thm}

\begin{remark}
                       \label{main remark} \rm

(i) Note that Theorem \ref{main thm 2} requires stronger conditions than those in Theorem \ref{main thm 1}; $\sigma^{ik}$ is assumed to be zero, and  the regularity condition of $\mu$ is   stronger  since $|\mu|_{B^{|\gamma+1|_+}}\leq |\mu|_{B^{|\gamma+1|^+}}$.

(ii) However Theorem \ref{main thm 2} gives better regularity results of solutions; let    $\gamma+2-d/p>0$ and $u$ be the solution in the above theorems. Then from the embedding $H^{\gamma+2}_p \subset C^{\gamma+2-d/p}$, it follows that
$$
\E \left[ \int^T_0 |u|^p_{C^{\gamma+2-d/p}}ds\right] \leq c
\left(\|f\|_{\bH^{\gamma}_p(T)}+\|g\|_{\bH^{\gamma+1+\varepsilon}_p(T,
\ell^2)}+\|u_0\|_{U^{\gamma+2}_p}\right).
$$
\end{remark}

\mysection{SPDEs with constant coefficients}

In this section we consider the equation
\begin{equation}      \label{eqn constant}
du= \left(a^{ij}u_{x^ix^j}+ f \right)\, dt + \left( \sigma^{ik}
u_{x^i}+g^{k} \right) \,dZ^{k}_{t},
\end{equation}
where the coefficients $a^{ij},\sigma^{ik}$ are independent of $x$.
 Recall that $\delta_0>0$ is the constant so that $(a^{ij})\geq \delta_0  I_{d\times d}$.

Let $T_t$ denote the semigroup associated with the Laplacian
$\Delta$ on $\bR^d$, that is,
$$
T_tf(x)=P_t*f(x), \qquad \hbox{where} \quad
 P_t(x)=(2\pi t)^{-d/2}e^{-|x|^2/(2t)}.
$$

\begin{lemma}
               \label{lemma 3.1}
Let $p\in [2,\infty)$ and $g=(g^1,g^2,\cdots)\in L^p((0,T)\times
\bR^d, \ell^2)$. Then
\begin{equation}
                       \label{e:12.1}
\int_{\bR^d}\int^T_0\left(\int^t_0 | DT_{t-s}g(s,x) |^2_{\ell^2} \,
ds\right)^{p/2}dtdx\leq c(d,p)\int_{\bR^d}\int^T_0|g(t,
x)|^p_{\ell^2}\, dtdx.
\end{equation}
\end{lemma}
\pf See Lemma 4.1 in \cite{Kr99}.

\qed

\begin{lemma}
                   \label{lemma 3.2}
Let $p\in (2,\infty)$ and $f\in L^p((0,T)\times \bR^d)$. Then for
any $\varepsilon>0$,
\begin{equation}
                       \label{e:12.1.1}
\int_{\bR^d}\int^T_0\int^t_0|DT_{t-s}f(s,x)|^p\,dsdt\,dx\leq
c\int^T_0\|f(t, \cdot)\|^p_{H^{\varepsilon}_p}\, dt,
\end{equation}
where $c=c(d,p,\varepsilon,T)$.
\end{lemma}

\pf Let $q>p$ be chosen so that $1/p=(1-\varepsilon)/2+\varepsilon
/q$, and define
 an operator $\cA$ by
$$
\cA f(t,s,x)=\begin{cases} D T_{t-s}f  \quad &\hbox{if } s<t, \\
 0 &\hbox{otherwise}. \end{cases}
 $$
 Then, due to Lemma \ref{lemma 3.1} and the inequality
 $\|T_{t-s}Df\|_q\leq \|Df\|_q$,
the linear mappings
 $$
\cA: L^2([0,T],L^2(\bR^d)) \to L^2([0,T]\times [0,T]\times \bR^d)
$$
 and
$$
\cA: L^q([0,T],H^1_q) \to L^q([0,T]\times [0,T]\times \bR^d)
$$
are bounded. Thus the lemma follows from the interpolation theory;
see, for instance, \cite[Theorem 5.1.2]{BL}.
\qed

\begin{remark} \rm
We  suspect  that (\ref{e:12.1.1})  is not true if $\varepsilon=0$,
and this is one of main reasons why we assumed $\sigma^{ik}=0$ in
Theorem \ref{main thm 2}.
\end{remark}

Here are the main results of this section.

\begin{thm}
                   \label{thm constant L_2}
For  every $f\in \bH^{\gamma}_2(T)$, $g\in \bH^{\gamma+1}_2(T,
\ell^2)$, $u_0\in U^{\gamma+2}_2$ and $T>0$, equation (\ref{eqn
constant}) with initial data $u_0$ has a unique solution $u\in
\cH^{\gamma+2}_2(T)$, and
\begin{equation}
                    \label{new 1}
\|u_{x}\|_{\bH^{\gamma+1}_2(T)}\leq
c\left(\|f\|_{\bH^{\gamma}_2(T)}+\|g\|_{\bH^{\gamma+1}_2(T,
\ell^2)}+\|u_0\|_{U^{\gamma+2}_2}\right)
\end{equation}
\begin{equation}
                         \label{e:2-2}
\|u\|_{\bH^{\gamma+2}_2(T)}\leq
ce^{cT}\left(\|f\|_{\bH^{\gamma}_2(T)}+\|g\|_{\bH^{\gamma+1}_2(T,
\ell^2)}+\|u_0\|_{U^{\gamma+2}_2}\right),
\end{equation}
where $c=c(\delta_0,d,K)$ is independent of $T$.
\end{thm}

\pf  Owing to Remark \ref{remark 11.22},  we may assume $\gamma=-1$.
Indeed, suppose that the theorem holds when $\gamma=-1$. Then it is
enough to notice that $u\in \cH^{\gamma+2}_2(T)$ is a solution of
the equation if and only if $v:=(1-\Delta)^{(\gamma+1)/2}u$ is a
solution of the equation with $\bar{f}:=(1-\Delta)^{(\gamma+1)/2}f$,
$\bar{g}:=(1-\Delta)^{(\gamma+1)/2}g$ and
 $\bar{u}_0:=(1-\Delta)^{(\gamma+1)/2}u_0$ in place of
$f,g$ and $u_0$, respectively, and
\begin{eqnarray*}
\|u\|_{\bH^{\gamma+2}_2(T)}=\|v\|_{\bH^1_2(T)} &\leq&
c\left(\|\bar{f}\|_{\bH^{-1}_2(T)} +\|\bar{g}\|_{\bH^{0}_2(T,
\ell^2)}
+\|\bar{u}_0\|_{U^{1}_2} \right) \\
&=&c \left( \|f\|_{\bH^{\gamma}_2(T)}+\|g\|_{\bH^{\gamma+1}_2(T,
\ell^2)}+\|u_0\|_{U^{\gamma+2}_2} \right).
\end{eqnarray*}

Since the coefficients $a^{ij}$ are independent of $x$, equation (\ref{eqn constant}) can be  rewritten as
$$
du= \left(\frac{\partial}{\partial x^i} \left( a^{ij}u_{x^j}\right)
 + f \right)\, dt + \left( \sigma^{ik} u_{x^i}+g^{k} \right) \,dZ^{k}_{t}.
$$
 By Remark 2.9 and Theorem 2.10 in \cite{C-K},
this equation has a unique solution $u\in \bH^1_2(T)$, and
furthermore  there is a constant $c>0$ independent of $T>0$ so that
$$
\|u_{x}\|_{\bL_2(T)}\leq c \left( \|f\|_{\bH^{-1}_2(T)}+
\|g\|_{\bH^0_2(T,\ell^2)}+ \|u_0\|_{U^{1}_2}\right).
$$
$$
   \|u\|_{\bH^1_2(T)}\leq  ce^{cT}
   \left( \|f\|_{\bH^{-1}_2(T)}+ \|g\|_{\bH^0_2(T,\ell^2)}+ \|u_0\|_{U^{1}_2}\right).
$$
The theorem is proved.
\qed

\begin{thm}
                   \label{thm constant L_p}
Let $p>2$, $\sigma^{ik}=0$ for each $i,k$,  $\varepsilon>0$ and
$T>0$. For  every $f\in \bH^{\gamma}_p(T)$, $g\in
\bH^{\gamma+1+\varepsilon}_p(T, \ell^2)$ and $u_0\in
U^{\gamma+2}_p$, equation (\ref{eqn constant}) with initial data
$u_0$ has a unique solution $u\in \cH^{\gamma+2}_p(T)$, and
\begin{equation}
                         \label{e:L_p}
\|u\|_{\cH^{\gamma+2}_p(T)}\leq
c(\|f\|_{\bH^{\gamma}_p(T)}+\|g\|_{\bH^{\gamma+1+\varepsilon}_p(T,
\ell^2)}+\|u_0\|_{U^{\gamma+2}_p}),
\end{equation}
where $c=c(\delta_0,d,p,\varepsilon,K,T)$.
\end{thm}

\pf Again by Remark \ref{remark 11.22}, we only need to prove the
theorem for $\gamma=-1$. Since the uniqueness of solution of equation (\ref{eqn constant}) follows from
the uniqueness result of deterministic  equations, we only
need to show
that there is a solution $u\in \cH^1_p(T)$ and $u$ satisfies
(\ref{e:L_p}).

{\bf Step 1}. First we prove the theorem for the stochastic heat
equation:
\begin{equation}
                        \label{e: heat eqn}
  du=\Delta u dt + \sum_{k=1}^{\infty} g^k dZ^k_t, \quad u(0)=0.
  \end{equation}
Using a standard approximation arguments, without loss of
generality, we may and do assume that $g^k=0$ for all  $k\geq N_0+1$
and that
$$
g^k(t,x)=\sum_{i=0}^{m(k)}I_{(\tau^k_i,\tau^k_{i+1}]}(t)g^{ki}(x),
$$
where $\tau^k_i$ are bounded stopping times and $g^{ki}(x)\in
 C^{\infty}_c (\bR^d)$. Define
$$
v(t,x):=\sum_{k=1}^{N_0} \int^t_0
 g^k(s, x) dZ^k_s
=\sum_{k=1}^{N_0}\sum_{i=1}^{m(k)} g^{ki}(x)(Z^k_{t\wedge
\tau^k_{i+1}}-Z^k_{t\wedge \tau^k_{i}})
$$
and
$$
u(t,x):=v(t,x)+\int^t_0 T_{t-s}\Delta v_s \, ds.
$$
 Then $d(u-v)=(\Delta(u-v)+\Delta v)dt=\Delta u dt$.
 Therefore
$$
du=\Delta u dt +dv=\Delta u dt +g^k dZ^k_t.
$$
Also by stochastic Fubini theorem,
almost surely,
\begin{eqnarray*}
u(t,x)&=&v(t,x)+\sum_{k=1}^{N_0}\int^t_0 \int^s_0T_{t-s}\Delta
g^k (r, x)dZ^k_r ds\\
&=&v(t,x)-\sum_{k=1}^{N_0}\int^t_0\int^t_r\frac{\partial}{\partial
s}T_{t-s}g^k (r, x) ds dZ^k_r\\
 &=&\sum_{k=1}^{N_0} \int^t_0 T_{t-s}g^k dZ^k_s.
\end{eqnarray*}
Similarly,
$$
\frac{\partial}{\partial x^i}u(t,x)=\sum_{k=1}^{N_0}\int^t_0
D_iT_{t-s}g^kdZ^k_s.
$$
Thus by
 Burkholder-Davis-Gundy's
inequality and (\ref{e: 12.6.4}), we have
\begin{eqnarray*}
\E \left[  |u_x(t,x)|^p \right]
&\leq& c\E\left[ \left(\sum_{k=1}^{N_0}\int^t_0\int
|DT_{t-s}g^k|^2|z|^2N^k(dz,ds)\right)^{p/2}\right]  \\
&\leq&  c(p)\, \E\left[ \left( \int^t_0\sum_{k=1}^{\infty}|DT_{t-s}\wh{g}^k|^2ds\right)^{p/2}
+ \int^t_0\sum_{k=1}^{\infty}|DT_{t-s}\wh{g}^k|^pds\right].
 \end{eqnarray*}
By Lemma \ref{lemma 3.1}, Lemma \ref{lemma 3.2} and the inequality
$\sum_k |a_k|^p\leq (\sum_k |a_n|^2)^{p/2}$,
\begin{equation}
                               \label{e:12.2.3}
\E \left[ \int^T_0\|Du\|^p \, dt \right]\leq
c\E \left[ \int^T_0\|g\|^p_{H^{\varepsilon}_p(\ell^2)} \, dt\right].
\end{equation}
 Next we prove (\ref{e:L_p}). As before,
\begin{equation}
                                 \label{eqn 8.6.5}
\E \left[ |u(t,x)|^p \right] \leq
c(p)\, \E\left[ \left( \int^t_0\sum_{k=1}^{N_0}|T_{t-s}\wh{g}^k|^2ds\right)^{p/2}
+ \int^t_0\sum_{k=1}^{N_0}|T_{t-s}\wh{g}^k|^p\, ds\right].
 \end{equation}
 Since $(\sum_{k=1}^{N_0}|a_n|^2)^{p/2}\leq
 N(N_0,p)\sum_{k=1}^{N_0}|a_n|^p$ and $\|T_tf\|_p\leq
 \|f\|_p$, it easily follows that $u\in \bH^0_p(T)$,
 and consequently $u\in \cH^1_p(T)$.
 By Theorem \ref{theorem banach} and (\ref{e:12.2.3}),
 $$
 \E \left[\sup_{s\leq t}\|u\|^p_{H^{-1}_p}\right]
 \leq c(d,p) \left(\|\Delta
 u\|^p_{\bH^{-1}_p(t)}+
  \|g\|^p_{\bH^{-1}_p(t,\ell^2)}\right)
 \leq c\|g\|^p_{\bH^{\varepsilon}_p(t,\ell^2)}.
 $$
This  together with \eqref{e:12.2.3}  and the inequality
 $$
 \|u\|_{H^1_p}=\|(1-\Delta)u\|_{H^{-1}_p}\leq  \|u\|_{H^{-1}_p}+
 \|\Delta u\|_{H^{-1}_p}\leq \|u\|_{H^{-1}_p}+
 \|Du\|_p
 $$
prove (\ref{e:L_p}). Before we move to next step, we emphasize that \eqref{eqn 8.6.5}
has nothing to do with \eqref{e:L_p},
and it is used only to show that $\|u\|_{\bH^0_p(T)}<\infty$.

{\bf{Step 2}}. General case. Let $v\in \cH^1_p(T)$ be the solution
of equation (\ref{e: heat eqn}), where the existence of the solution is obtained in Step 1.
Also let $\bar{u}\in \cH^1_p(T)$ be the
solution of the following equation (see \cite[Theorem 4.10]{Kr99})
$$
d\bar{u}=(a^{ij}\bar{u}_{x^ix^j}+f+a^{ij}v_{x^ix^j}-\Delta v)dt,
\quad \bar{u}(0)=u_0.
$$
Then by Step 1 and \cite[Theorem 4.10]{Kr99},
$$
\|v\|_{\bH^{1}_p(T)}\leq
c\|g\|_{\bH^{\varepsilon}_p(\ell^2)},
$$
$$
\|\bar{u}\|_{\bH^{1}_p(T)}\leq
c(\|v_{xx}\|_{\bH^{-1}_p(T)}+\|f\|_{\bH^{-1}_p(T)}+\|u_0\|_{U^{1}_p}).
$$
Note that $u:=\bar{u}+v$ satisfies
$$
du=(a^{ij}u_{x^ix^j}+f)dt +g^kdZ^k_t, \quad u(0)=u_0
$$
and estimate (\ref{e:L_p}) follows.
\qed

\mysection{Proof of Theorem \ref{main thm 1} and Theorem \ref{main thm 2}}

First we prove the following lemmas.

\begin{lemma}
                               \label{lemma 8.1}
For $c>0$, let   $Z^k_{t}(c):=c^{-1}Z^k_{c^2t}$ and $\nu^k_c$ be  the L\'e{}vy measure of $Z^k_{t}(c)$. Then
$$
\int_{\bR} z^2 \nu^k(dz) =\int_{\bR} z^2 \nu^k_c(dz).
$$
\end{lemma}

\pf Denote $N^k_c(t,A)=\#\{s\leq t; \Delta Z^k_{s}(c)\in A \}$. Then $N^k_c(t,A)=N^k(c^{2}t,cA)$ and
$$
\nu^k_c(A)=\E \left[ N^k(c^2,cA)\right]=c^2 \nu^k(cA).
$$
Hence the lemma follows from a change of variables. Indeed, let $f(z)=cz$ and $h(z)=z^2/{c^2}$. Then $\nu^k_c \circ f^{-1}=c^2 \nu^k$, and so
$$
\int z^2 \nu^k_c(dz)=\int h(f(z)) \nu^k_c(dz)=\int h(z) c^2 \nu^k(dz)=\int z^2 \nu^k(dz).
$$
\qed

Consider equation (\ref{eqn constant}) with $Z^k_t(c)$ in place of
$Z^k_t$. It follows from \eqref{e:2.3}  with $p=2$ and Lemma
\ref{lemma 8.1} that Theorem \ref{thm constant L_2} and
(\ref{e:2-2}) holds for this new equation with the same constant
$C$. Recall the definition of $B^r$ from \eqref{e:2.11} and that
$r^+:=r+\kappa_0$.

\begin{lemma}
                                  \label{lemma 8.3}
 For $d\geq 1$, $p\geq 2$ and $\gamma \in \bR$, there
is a constant $N=N(d,p,\gamma)>0$ so that
for every $a\in B^{|\gamma|_+}$ and $u\in H^{\gamma}_p$,
$$
\|au\|_{H^{\gamma}_{p}}\leq N |a|_{B^{|\gamma|_+}}\,
\|u\|_{H^{\gamma}_p}.
$$
The same is true for $\ell^2$-valued functions $a$  in
$B^{|\gamma|_+}$.
\end{lemma}

\pf See Lemma 5.2 in \cite{Kr99}. \qed

\begin{lemma}
                               \label{lemma 8.2}
Let  $b^i=c= \mu^k=0$ and  suppose that there is a constant $L>0$ so
that
$$
|a^{ij}(t,\cdot) |_{B^{|\gamma|_+}} + +
|\sigma^i(t,\cdot)|_{B^{|\gamma+1|_+}}
\leq L.
$$
 Define
$$
\beta:=\sup_{\omega, x, y} \left( |a^{ij}(t,x)-a^{ij}(t,y)|+
|\sigma^{ik}(t,x)-\sigma^{ik}(t,y)|_{\ell^2} \right).
$$
Then there exists $\beta_0=\beta_0(d,\delta,K)>0$, {\bf{independent
of $L$}}, so that if $\beta \leq \beta_0$ then for  any solution of
$u\in \cH^{\gamma+2}_2(T)$ of equation \eqref{e:2.1} we have
\begin{equation}
                 \label{eqn lemma 8.2}
\|u\|_{\bH^{\gamma+2}_{2}(T)}\leq ce^{cT}\left(
\|f\|_{\bH^{\gamma}_2(T)}+\|g\|_{\bH^{\gamma+1}_2(T,\ell^2)}+
\|u_0\|_{U^{\gamma+2}_2}\right),
\end{equation}
 where $c=c(d,\delta,K,L)$.
\end{lemma}

\pf  Let $u\in \cH^{\gamma+2}_2(T)$ be a solution to equation
\eqref{e:2.1}. Denote
$$
a^{ij}_0(t)=a^{ij}(t,0), \quad \sigma^{ik}_0(t)=\sigma^{ik}(t,0),
$$
$$
f_0=(a^{ij}-a^{ij}_0)u_{x^ix^j}+f, \quad g^k_0=(\sigma^{ik}-\sigma^{ik}_0)u_{x^i}+g^k,
$$
$$
C_0=\sup_{\omega,t}\left(|a^{ij}-a^{ij}_0|_{B^{|\gamma|_+}}+|\sigma^{i}-\sigma^i_0|_{B^{|\gamma+1|_+}}\right).
$$
Then  $\displaystyle du=(a^{ij}_0u_{x^ix^j}+f_0) dt +
(\sigma^{ij}_0u_{x^i}+g^k_0)dZ^k_t$. So
 by Theorem \ref{thm constant L_2},
$$
\|u_{x}\|_{\bH^{\gamma+1}_2(T)}\leq
c(d,\delta,K)\left(\|f_0\|_{\bH^{\gamma}_2(T)}+\|g_0\|_{\bH^{\gamma+1}_2(T,\ell^2)}+
\|u_0\|_{U^{\gamma+2}_2}\right).
$$
By Lemma \ref{lemma 8.3}
$$
\|(a^{ij}-a^{ij}_0)u_{x^ix^j}\|_{H^{\gamma}_2}\leq N(d,\gamma)
|a^{ij}-a^{ij}_0|_{B^{|\gamma|_+}}\|u_{x^ix^j}\|_{H^{\gamma}_2}\leq
NC_0\|u_{xx}\|_{H^{\gamma}_2}\leq NC_0\|u_x\|_{H^{\gamma+1}_2},
$$
and similarly
$$
\|(\sigma^i-\sigma^i_0)u_{x^i}\|_{H^{\gamma+1}_2(\ell^2)}\leq
NC_0\|u_x\|_{H^{\gamma+1}_2}.
$$
Thus the lemma follows if $cNC_0\leq 1/4$.  For $m\geq 1$,  denote
$a^{ij}_m(t,x):=a^{ij}(t/{m^2},x/m)$ and
$\sigma^{ik}_m(t,x):=\sigma^{ik}(t/{m^2},x/m)$.  Then we have
$$
|a^{ij}_m(t,\cdot)-a^{ij}_m(t,0)|_{B^{|\gamma|_+}} \leq
\beta+m^{-(|\gamma|_+ \wedge 1)}C_0,
$$
and we can drop the second term on the right if $\gamma=0$. Also we
have a similar inequality for $\sigma^{ik}_m$. Observe that
 $u_m(t,x):=u(t/{m^2},x/m)$ satisfies
$$
du_m=(a^{ij}_m (u_m)_{x^ix^j}+f_m)\,dt
+(\sigma^{ik}_m (u_m)_{x^i}+g^{k}_m)\,dZ^k_t(m^{-1}),
$$
where $f_m(t,x):=m^{-2}f(t/{m^2},x/m)$ and
$g^{ik}_m(t,x):=m^{-1}g^{ik}(t/{m^2},x/m)$. Then it follows from the
above calculations and  Lemma \ref{lemma 8.1} that for $\beta$
sufficiently small and $m$ sufficiently large,
$$
\|u_{mx}\|_{\bH^{\gamma+1}_2(mt)}\leq
c\left(\|f_m\|_{\bH^{\gamma}_2(mt)}
+\|g_m\|_{\bH^{\gamma+1}_2(mt,\ell_2)}+
\|u_0(\cdot/{m})\|_{U^{\gamma+2}_2}\right)
$$
for each $t\leq T$. Also, since $\|\cdot\|_{H^{\gamma}_p}$ norms of
$u(t/{m^2},x/m)$ and $u(t,x)$ are comparable, one gets inequality
(\ref{new 1}) for each $t\leq T$ in place of $T$.

By (\ref{new 1}) and the inequality
$$
\|u\|_{H^{\gamma+2}_2}=\|(1-\Delta)u\|_{H^{\gamma}_2}\leq
\|u\|_{H^{\gamma}_2}+c\|u_x\|_{H^{\gamma+1}_2}
$$
it follows that for each $t\leq T$,
$$
\|u\|^2_{\cH^{\gamma+2}_{2}(t)}\leq
c\left(\|u\|^2_{\bH^{\gamma}_2(t)}+\|f\|^2_{\bH^{\gamma}_2(t)}+\|g\|^2_{\bH^{\gamma+1}_2(t)}
+\|u_0\|^2_{U^{\gamma+2}_2}\right).
$$
This, (\ref{new 2}) and Gronwall's inequality prove (\ref{eqn lemma
8.2}).

\qed

\begin{lemma}
                            \label{lemma 8.4}
Let $\zeta_n\in C^{\infty}, n=1,2,\cdots$. Assume that for any
multi-index $\alpha$,
\begin{equation}
                          \label{eqn 03.30.1}
\sup_{x}\sum_n |D^{\alpha}\zeta_n(x)|\leq M(\alpha),
\end{equation}
where $M(\alpha)$ are some constants. Then there exists a constant
$N=N(d,n, \gamma, p, M)$ such that for any $f\in H^{\gamma}_p$,
$$
\sum_n \|\zeta_nf\|^p_{H^{\gamma}_p}\leq N\|f\|^p_{H^{\gamma}_p}.
$$
If, in addition,
\begin{equation}
                 \label{eqn 03.30.2}
\sum_n|\zeta_n(x)|^p >c>0,
\end{equation} then
$$
\|f\|^p_{H^{\gamma}_p}\leq  N(d,n, \gamma, p, M,c)
 \sum_n \|\zeta_nf\|^p_{H^{\gamma}_p}.
$$
 \end{lemma}
\pf See Lemma 6.7 in \cite{Kr99}.
\qed

\medskip

\noindent {\bf{Proof of Theorem \ref{main thm 1}}} \ \
 Let $\beta_0>0$ be the constant in Lemma \ref{lemma 8.2}.
In view of Theorem \ref{thm constant L_2} and the method of
continuity (see the proof of \cite[Theorem 2.11]{C-K}), we only need
to show that a priori estimate \eqref{eqn main L_2} holds given that
a solution $u\in \cH^{\gamma+2}_2(T)$ already exists. Take $\beta_0$
from Lemma \ref{lemma 8.2}. Since $a^{ij}, \sigma^i$ are uniformly
continuous we can fix $\delta_0$ so that
$$
|a^{ij}(t,x)-a^{ij}(t,y)|+
|\sigma^{i}(t,x)-\sigma^i(t,y)|_{\ell^2}\leq \beta_0/4
$$
if $|x-y|\leq 2\delta_0$. Fix   a smooth function $\zeta\in
C^{\infty}_c(B_1(0))$ so that $0\leq \zeta\leq 1$ and $\zeta(x)=1$
if $|x|\leq 1/2$. Take a sequence of smooth functions
$\{\zeta_n:n=1,2,\cdots\}$ so that
$\zeta_n=\zeta(\frac{2(x-x_n)}{\delta_0})$ for some $x_n\in \bR^d$,
and (\ref{eqn 03.30.1}) and (\ref{eqn 03.30.2}) hold.

 Then by Lemma
\ref{lemma 8.4}, for each $t\leq T$,
\begin{equation}
                                      \label{eqn 8.3.1}
\|u\|^2_{\bH^{\gamma+2}_2(t)}\leq N\sum_n \|u\zeta_n\|^2_{\bH^{\gamma+2}_2(t)}.
\end{equation}
Denote $\xi_n(x)=\zeta(\frac{x-x_n}{\delta_0})$ and
$$
a^{ij}_n(t,x)=\xi_n(x)a^{ij}(t,x)+(1-\xi_n (x))a^{ij}(t,x_n), \quad \sigma^{ik}_n(t,x)=\xi_n(x)\sigma^{ik}(t,x)+(1-\xi_n (x))\sigma^{ik}(t,x_n).
$$
Then $a_n$ and $\sigma_n$ satisfy (\ref{eqn 8.5.1})  with the same constants $\delta,K$,
$$
|a^{ij}_n(t,x)-a^{ij}_n(t,y)|+
|\sigma^{i}_n(t,x)-\sigma^i_n(t,y)|_{\ell^2}\leq \beta_0, \quad
\forall \omega,t,x,y
$$
and $u\zeta_n$ satisfies
$$
d(u\zeta_n)=(a^{ij}_n(u\zeta_n)_{x^ix^j}+f_n)dt+(\sigma^{ik}_n(u\zeta_n)_{x^i}+g^k_n)dZ^k_t,
$$
where
$$
f_n=-2a^{ij}u_{x^i}\zeta_{nx^j}-a^{ij}u\zeta_{nx^ix^j}+b^iu_{x^i}\zeta_n+cu\zeta_n+f\zeta_n
$$
$$
g^k_n=-\sigma^{ik}u\zeta_{nx^i}+\mu^ku\zeta_n+g^k\zeta_n.
$$
By Lemma \ref{lemma 8.2} and Lemma \ref{lemma 8.3}
\begin{eqnarray*}
&& \|u\zeta_n\|^2_{\bH^{\gamma+2}_2(t)}
 \leq  c \left( \|f_n\|^2_{\bH^{\gamma}_2(t)}
 + \|g_n\|^2_{\bH^{\gamma+1}_2(t,\ell^2)}+ \|u_0\zeta_n\|^2_{U^{\gamma+2}_2} \right) \\
&\leq&  c \left( \|u_x\zeta_{nx}\|^2_{\bH^{\gamma}_2(t)}+ \|u\zeta_{nxx}\|^2_{\bH^{\gamma}_2(t)}+c\|u\zeta_{nx}\|^2_{\bH^{\gamma+1}_2(t)}+
 \|f\zeta_n\|^2_{\bH^{\gamma}_2(t)}+
 \|g\zeta_n\|^2_{\bH^{\gamma+1}_2(t,\ell^2)}+ \|u_0\zeta_n\|^2_{U^{\gamma+2}_2}\right).
\end{eqnarray*}
Thus by (\ref{eqn 8.3.1}) and Lemma \ref{lemma 8.4}
\begin{equation}\label{e:4.3}
\|u\|^2_{\bH^{\gamma+2}_2(t)}\leq c \left(
\|u\|^2_{\bH^{\gamma+1}_2(t)} + \|f\|^2_{\bH^{\gamma}_2(t)}+
\|g\|^2_{\bH^{\gamma+1}_2(t,\ell^2)}+
\|u_0\|^2_{U^{\gamma+2}_2}\right).
\end{equation}
By definition of the space $\cH^{\gamma+2}_2(t)$ and  Lemma \ref{lemma 8.3}
\begin{eqnarray*}
\|u\|_{\cH^{\gamma+2}_2(t)}&:=& \|u\|_{\bH^{\gamma+2}_2(t)}
+\|a^{ij}u_{x^ix^j}+b^iu_{x^i}+cu+f\|_{\bH^{\gamma}_2(t)}
+\|\sigma^{i}u_{x^i}+\mu u+g\|_{\bH^{\gamma+1}_2(t,\ell^2)}+\|u_0\|_{U^{\gamma+2}_2} \\
&\leq& N\|u\|_{\bH^{\gamma+2}_2(t)}+\|f\|_{\bH^{\gamma}_2(t)}
+\|g\|_{\bH^{\gamma+1}_2(t,\ell^2)}+\|u_0\|_{U^{\gamma+2}_2}.
\end{eqnarray*}
 This together with \eqref{e:4.3}, the embedding inequality (see,
 for instance, \cite{T})
\begin{equation}
            \label{embedding}
\|u\|^2_{H^{\gamma+1+\beta}_2}\leq \varepsilon
\|u\|^2_{H^{\gamma+2}_2}+c(\varepsilon,\beta)\|u\|^2_{H^{\gamma}_2},
\quad \forall \beta\in [0,1)
\end{equation}
and Theorem \ref{theorem banach} yields that
$$
\|u\|^2_{\cH^{\gamma+2}_2(t)}\leq
c\int^t_0\|u\|^2_{\cH^{\gamma+2}_2(s)}ds
+c\left(\|f\|^2_{\bH^{\gamma}_2(T)}+\|g\|^2_{\bH^{\gamma+1}_2(T,\ell^2)}+
\|u_0\|^2_{U^{\gamma+2}_2}\right).
$$
The a priori estimate now follows from Gronwall's inequality. The theorem is proved.  \qed

\medskip

\noindent {\bf{Proof of Theorem \ref{main thm 2}}}\
 Again,  in view of Theorem \ref{thm constant L_p} and the method of continuity,  we only need to
prove that a priori estimate (\ref{eqn main L_p}) holds given that a solution $u\in \cH^{\gamma+2}_p(T)$ already exists. Without loss of generality we assume that
$\varepsilon< (\kappa_0\wedge 1)$ since once (\ref{eqn main L_p}) holds for some small $\varepsilon$ then it holds for any $\varepsilon'\geq \varepsilon$. One can prove the theorem by modifying the proof of Theorem \ref{main thm 1}. But since $\sigma^{ik}$ is assumed to be zero, the proof of this theorem is much easier.

Let $v\in \cH^{\gamma+2}_p(T)$ be the solution of
$$
dv=\Delta v dt + ( \mu^k u+g^k)dZ^k_t, \quad v(0)=u_0.
$$
The existence of the solution of the above equation is guaranteed by Theorem \ref{thm constant L_p}. By Theorem \ref{thm constant L_p} and Lemma \ref{lemma 8.3},
\begin{equation}\label{e:4.4}
\|v\|_{\bH^{\gamma+2}_p(t)}\leq c \left(\| \mu
u+g\|_{\bH^{\gamma+1+\varepsilon}_p(t,\ell^2)} +
\|u_0\|_{U^{\gamma+2}_p} \right)\leq c \left(
\|u\|_{\bH^{\gamma+1+\varepsilon}_p(t)}+
\|g\|_{\bH^{\gamma+1+\varepsilon}_p}+
\|u_0\|_{U^{\gamma+2}_p}\right).
\end{equation}
Note that $\bar{u}:=u-v\in \cH^{\gamma+2}_p(T)$ satisfies
\begin{equation}\label{e:4.5}
d\bar{u}=(a^{ij}\bar{u}_{x^ix^j}+b^i\bar{u}_{x^i}+c\bar{u}+\bar{f})dt, \quad \bar{u}(0)=0,
\end{equation}
where $\bar{f}=a^{ij}v_{x^ix^j}-\Delta v+b^iv_{x^i}+cv+f$. By Theorem 5.2 in \cite{Kr99},
$$
\|\bar{u}\|_{\bH^{\gamma+2}_p(t)}\leq c\|\bar{f}\|_{\bH^{\gamma}_p(t)}\leq c\left( \|v\|_{\bH^{\gamma+2}_p(t)}+ \|f\|_{\bH^{\gamma}_p(t)}\right).
$$
Consequently, for each $t\leq T$,  by \eqref{e:4.4} and
\eqref{e:4.5}
$$
\|u\|_{\bH^{\gamma+2}_p(t)}
\leq c \left( \| \bar u\|_{\bH^{\gamma+2}_p(t)}
 + \| v\|_{\bH^{\gamma+2}_p(t)}\right)
\leq c \left( \|u\|^p_{\bH^{\gamma+1+ \varepsilon}_p(t)}+
\|f\|^p_{\bH^{\gamma}_p(t)} +
\|g\|^p_{\bH^{\gamma+1+\varepsilon}_p(t)}+
\|u_0\|^p_{U^{\gamma+2}_p} \right).
$$
This together with the embedding inequality (see (\ref{embedding}))
$$
\|u\|_{H^{\gamma+1+\varepsilon}_p}\leq \delta \|u\|_{H^{\gamma+2}_p}+c(\delta,\varepsilon)\|u\|_{H^{\gamma}_p},
$$
yields that
$$
\|u\|_{\bH^{\gamma+2}_p(t)}\leq c \left( \|u\|^p_{\bH^{\gamma}_p(t)}
+\|f\|^p_{\bH^{\gamma}_p(t)}+
\|g\|^p_{\bH^{\gamma+1+\varepsilon}_p(t)} +
\|u_0\|^p_{U^{\gamma+2}_p}\right).
$$
As in the proof of Theorem \ref{main thm 1}, this easily leads to (\ref{eqn main L_p}). The theorem is proved. \qed

\end{document}